\documentclass[leqno]{amsart}
\usepackage{amsmath,amsthm,amscd}

\def\a{\alpha}
\def\b{\beta}
\def\d{\delta}
\def\g{\gamma}
\def\s{\sigma}

\def\supth{{\text{th}}}
\def\f{\frac}
\def\op{\oplus}
\def\ov{\overline}
\def\wt{\widetilde}
\def\wh{\widehat}
\def\lra{\longrightarrow}
\def\SC{\mathcal C}
\def\SD{\mathcal D}
\def\SE{\mathcal E}
\def\SF{\mathcal F}
\def\SO{\mathcal O}
\def\SS{\mathcal S}
\def\ST{\mathcal T}
\def\SU{\mathcal U}
\def\SZ{\mathcal Z}
\newcommand{\BC}{\ensuremath{\mathbf C}}
\newcommand{\BD}{\ensuremath{\mathbf D}}
\newcommand{\BH}{\ensuremath{\mathbf H}}
\newcommand{\BP}{\ensuremath{\mathbf P}}
\newcommand{\BR}{\ensuremath{\mathbf R}}
\newcommand{\BZ}{\ensuremath{\mathbf Z}}
\def\oper{\operatorname}
\newcommand{\Div}{\operatorname{Div}}

\def\ger{\mathfrak g}
\def\sd#1{\SD{'{^#1}}}
\def\sdk{\SD{'{^{k-1}}}}
\def\spark{\oper{spark}}

\def\tor{\oper{tor}}
\def\free{\oper{free}}

\def\hypers{\oper{hyperspark}}
\def\sm{\oper{smooth}}
\def\loc{L^1_{\text{loc}}}
\def\hom{\oper{Hom}}
\def\CDR{\widetilde{C}_{\oper{deR}}}

\def\dbar{\overline{\partial}}

\numberwithin{equation}{section}

\newtheorem{thm}{Theorem}[section]
\newtheorem{prop}[thm]{Proposition}
\newtheorem{lem}[thm]{Lemma}
\newtheorem*{haf}{(HAF)}

\theoremstyle{definition}
\newtheorem{defn}[thm]{Definition}
\newtheorem*{note}{Note}

\begin{document}

\title{From Sparks to Grundles --- Differential Characters}
\author{Reese Harvey and Blaine Lawson}
\date{December 1, 2003}
\maketitle

\centerline{\bf Abstract} \medskip
  \font\abstractfont=cmr10 at 10 pt

{{\parindent=.9in\narrower\abstractfont \noindent
We introduce a new homological machine for the study of secondary
geometric invariants. The objects, called spark complexes, occur in many
areas of mathematics.  The theory is applied here to establish the
equivalence of a large family of spark complexes which appear naturally in
geometry, topology and physics.  These complexes are quite different. Some
of them are purely analytic, some are  simplicial, some are of  \v
Cech-type, and many are mixtures.  However, the  associated theories of
secondary invariants are all shown to be canonically isomorphic.
Numerous applications and examples are explored.

}}

\vskip .4in

\section*{Contents}

\quad1. Homological Sparks

\quad2. Sparks (DeRham-Federer)

\quad3. Hypersparks

\quad4. Smooth Hypersparks

\quad5. Grundles

\quad6. Cheeger-Simons Sparks

\quad7. Differential Characters

\quad8. Holonomy Maps

\quad9. Grundle Holonomy

\quad10. Further Spark Complexes

\quad11. Applications and Examples


\newpage
\section*{Introduction}

The point of this paper is to present a general homological apparatus
for the study of sparks, differential characters, gerbes and analogous
objects.  This overarching theory encompasses all the known
theories equivalent to differential characters, as well as many new ones,
and it establishes this equivalence.  The basic algebra is elementary and
self-contained.  Nevertheless, it enables one to rigorously establish
canonical isomorphisms between the quite disparate approaches to the
theory: the original Cheeger-Simons formulation, the deRham-Federer theory
using currents, the \v Cech-deRham formulation connected with gerbes, and
many others.  It also provides a homological framework for the development
of new secondary theories which will be persued in   separate papers.

In the first section we introduce the concept of a  {\bf homological
spark complex} and its associated {\bf group of  homological spark
classes}. (These play roles analogous to cochain complexes and their
associated groups of cohomology classes in standard homological
algebra.)\   We establish a basic $3\times 3$ grid of short exact
sequences with the group of homological spark classes in the center.
The spark complexes of interest in this paper will be referred to
collectively as {\bf $(\BR,\BZ)$-spark complexes}. For these
complexes, the grid replicates the standard short exact sequences in
that theory.

We then introduce the notion of {\bf  compatible spark  complexes}
and show that compatible complexes have naturally isomorphic groups
of spark classes.  This is the tool that enables us to rigorously
unify the many approaches to the  $(\BR,\BZ)$ theory.

The next part of the paper is devoted to studying specific examples.
We shall present:
(HM) holonomy maps,
 (DC)  differential characters,
 (DFS) de-Rham-Federer spark classes,
 (HS) hyperspark classes,
 (SHS) smooth hyperspark classes,
 (GR) grundles (n-gerbes with connection),
 (CHS) cochain hyperspark classes,
 (CSCS)  Cheeger-Simons cochain spark classes,
(CCS) current cochain spark classes,
and prove that all of these concepts are equivalent.
The diagram 
$$
\begin{matrix}
\ &\ & \boxed{HS} &\  &\  &\boxed{CHS} &\ &\ \qquad \boxed{DC}\\
\  & \nearrow &\  & \nwarrow \ \ &  \nearrow &\  &  \nwarrow &\  
\nearrow    &\  \  \nwarrow\\
\boxed {DFS}  &\ &\   &\boxed{SHS} &\ &\   &\boxed{CSCS} &\ &\  \
\boxed{HM}\\
\  &\ &\   & \downarrow &&&\hskip40pt\ \nwarrow \ &\ &\nearrow \\
\  &\ &\  &\  \boxed{Grundles}\ &\  &\ &\ \ \ \ \ \ \ \  &\boxed{CCS}
\end{matrix}
$$
$$
\text{$(\BR,\BZ)$-Spark Complexes}
$$
provides a convenient summary. The arrows in this diagram
represent the direction of the inclusion maps at the level of
spark complexes. The paper is organized as follows.

In \S 2 we recall the deRham-Federer spark complexes  which first
appeared in [GS] and [Harris] and were systematically studied in
[HLZ].  There are  in fact many such complexes, from the most general
one involving all currents, to quite restricted ones in which only
forms with $\loc$-coefficients  appear as sparks [HLZ, \S 2].

In \S 3 we define the complex of {\sl hypersparks} whose compatibility
with all the deRham-Federer spark complexes establishes the isomorphisms
between their groups of spark classes.

In \S 4 we introduce {\sl smooth hypersparks} and show that they form a
compatible subcomplex of the hypersparks.

In \S 5 we present the concept of  {\sl grundles } of degree $k$ and the
notion of  {\sl gauge equivalence}. Grundles of degree 1 are just  \v
Cech presentations of hermitian line bundles with unitary connection, and
``gauge equivalence'' corresponds to the standard equivalence as unitary
bundles with connection.  Grundles of degree 2 similarly correspond to
(abelian) gerbes with connection, and ``gauge equivalence'' corresponds  to
the standard equivalence found in the literature [Hi], [Br].

Grundles of degree $k$ and smooth hypersparks of degree $k$ are essentially
the same.  they differ only in the lowest order component where one is
obtained from the other by choosing logarithms ( or inversely by
exponentiation). In particular we show that gauge equivalence classes of
grundles are isomorphic to the corresponding classes of smooth hypersparks
(which are isomorphic to groups of spark classes of all the aforementioned
spark complexes).

In \S 6 we consider {\sl Cheeger-Simons sparks} which are the sparks most
closely related to differential characters as originally defined in
[CS]. These sparks are shown to be compatible with the smooth hypersparks
via a larger complex which encompasses them both (cochain
hypersparks).

In \S 7 we examine Cheeger-Simons theory and relate it to all of the above.
Recall  that in  [CS]   a {\sl differential character} of
degree $k$ is defined to be a homomorphism
$$
h:Z_k(X)\ \lra\ \BR/\BZ\qquad \text{ such that }\ \ \ \  dh \equiv \phi
\mod \BZ
$$
for some smooth $(k+1)$-form $\phi$ on $X$. Here $Z_k(X)=\{c\in C_k(X):
\partial c=0\}$  where $C_k(X)$ is  the group of $C^{\infty}$ singular
$k$-chains in $X$ with integer coefficients, and the congruence  above
means that $h(\partial c)-\int_c\phi \in \BZ$ for all $c\in C_{k+1}(X)$.
We show that {\sl the group of spark classes for each  spark complex
appearing in \S\S 2--6 is functorially isomorphic  to the
space} DiffChar$^k(X)$.

>From a geometric point of view perhaps the most natural starting point for
the $(\BR,\BZ)$-theory is to define ``holonomy maps''.  Quite surprizingly
this concept does not seem to be part of the literature.  Note that every
$k$-chain $c \in C_k(X)$ defines a {\sl current}, that is, a linear
functional $\wt c$ on the space $\SE^k(X)$ of smooth $k$-forms on $X$, by
setting $\wt c (\psi)=\int_c \psi$. Let ${\wt Z}_k(X)$ be the group of
currents coming from cycles in $Z_k(X)$ in this way. (These  are {\sl
current  chains} in the sense of de Rham.) There is a short exact
sequence 
$$
0\ \to\ N_k(X)\ \to\ {Z}_k(X)\overset{\rho}{\lra}\ {\wt Z}_k(X)\ \to\
0
$$
where $\rho(z)=\wt z$ and $N_k(X)$ is the group of {\sl null cycles}
on $X$, the cycles over which all forms integrate to be 0.  Note for
example that a compact oriented submanifold $M\subset X$ of dimension $k$
defines an element $[M]\in {\wt Z}_k(X)$.  Among the many cycles $z$ with
$\rho(z)=[M]$ are those arising from  smooth triangulations of $M$.

In \S 8 we define the set of {\bf degree $k$ holonomy maps} to be
$$
{\bold H}{\bold o}{\bold l}(X)\ \equiv\
\{H\in\hom({\wt Z}_k(X), \BR/\BZ)\ :\
dH \equiv \phi \mod \BZ\ \ \text{for some } \phi \in \SE^{k+1}(X)\}
$$
In Theorem 8.2 we show
that {\sl every differential character $h$ as above can be written  in the
form $h=H\circ \rho$ for a homomorphism $H : {\wt Z}_k(X)\ \to\
\BR/\BZ$.} In other words there is an isomorphism
$$
\text{DiffChar}^k(X)\ \cong\ {\bold H}{\bold o}{\bold l}(X).
$$

It is well known that in degree 1 differential characters are isomorphic
to gauge equivalence classes of principal $S^1$-bundles with connection,
and that the isomorphism is given by the classical holonomy of the
connection. The above results generalize this statement to grundles of all
degrees.

The holonomy of any grundle over a compact oriented submanifold $M$ is
independent  of any choice of triangulation of $M$.  However, given a
triangulation it is interesting to find combinatorial formulas for the
holonomy in terms of the grundle data.  A general formula of this type
is derived in \S 9.  In degree 1 this formula computes the holonomy
around a loop by integration of the connection 1-forms on sub-arcs where
the bundle is trivialized and the change of trivialization at the endpoints
of the arcs.  The corresponding formula (9.5) for the holonomy of a gerbe
with connection involves integration of the gerbe data over the 0-,1- and
2-faces. For a detailed discussion of holonomy in this case see [MP].

There are many $(\BR,\BZ)$ homological spark complexes which were not
mentioned in  \S\S 2--6.  We briefly discuss some   useful ones in
Section 10.

We note that the discussion  here was restricted to the case of  ``integer
coefficients''. However,  one can easily replace $\BZ$ with any subring
$\Lambda\subset \BR$.  For example, in the discussion of \S 3, one
replaces  $C^*(\SU, \BZ) \subset C^*(\SU,\BR)$ with $C^*(\SU, \Lambda)
\subset C^*(\SU, \BR)$, and everything goes through.

Historically, the  spark equation and spark equivalence first appeared in
[GS] and [H] and were systematically developed in [HLZ], where
Poincar\'e-Pontrjagin duality for differential characters was established.
The concept of  smooth hypersparks can be found in the work of Freed-Witten
[FW], and Picken [P], and the concept of grundles (sometimes called {\sl
$n$-gerbes with connection}) can be found in many places, e.g. [Hi], [B].
Cheeger-Simons cochain sparks  appear explicitly in the papers of Cheeger
and Simons [C], [CS]  as ``liftings'' of differential characters
$\phi:Z_k(X)\to \BR/\BZ$ to homomorphisms $\wt\phi:Z_k(X)\to\BR$. However
they do not systematically study  the spark equation and spark equivalence.

In \S 11 we briefly discuss many applications and examples.  In each case
our main concern is to examine which points of view (see diagram above) are
most illuminating.

In this paper our attention is focused on $(\BR,\BZ)$-spark complexes.
However, there exist many spark complexes whose associated class groups
are  quite different and interesting.

For example on complex manifolds one can define   a {\bf $\dbar$-analogue
of differential characters } by applying (a slight generalization of) the
homological machinery developed here to spark  complexes involving forms
and currents of type $(0,q)$. There are several compatible formulations
of the theory which yield  the same $\dbar$-spark classes. One
formulation uses a smooth \v Cech - Dolbeault complex which is related
(as in this paper) to $\dbar$-grundles.  In low degrees these grundles
have beautiful geometric interpretations as in the classical case.
However, there is a quite different Dolbeault-Federer approach to
$\dbar$-spark classes which is purely analytic.  It enables us to
establish functoriality and to define a $*$-product which makes these
groups into a graded ring. The equivalence  of these spark complexes is
proved by a larger hyperspark complex compatible with them both, as in
this paper. One interesting by-product of the discussion is an acyclic
resolution of the sheaf $\SO^*$ by forms and currents.  All these
results will appear in [HL$_5$].

We emphasize that there exist many further interesting spark complexes.
Some involve $(p,q)$-forms and are related to Deligne cohomology.  Others
arise in a quite general way from interesting double complexes.   These
will be explored in yet another article [HL$_6$].

\begin{note} In discussing certain spark complexes we talk freely about
sheaves of germs of currents of various types.  For a unified approach to
the cohomology of current complexes (such as flat, integrally flat,
integral, etc.) the reader is referred to the Appendix A in [HZ$_1$].
\end{note}

\section{Homological Sparks} A \textbf{spark complex} consists of a
differential complex $F^*$, $d$, and two subcomplexes $E^*\subset F^*$ and
$I^*\subset F^*$ which satisfy the following properties:

\medskip
(A) $H^k(E^*)\cong H^k(F^*)$, $k=0$, $1,\dots$

(B) $E^k\cap I^k=\{0\}$, $k=1$, $2,\dots$\newline

\medskip
It follows that 
\begin{equation}
E^0\cap I^0=H^0(I^*)\subset H^0(E^*),\tag*{\qquad (C)}
\end{equation}
since $a\in E^0\cap I^0$ implies $da\in E^1\cap I^1=\{0\}$ and
$H^0(I^*)\subset H^0(F^*)=H^0(E^*)$.

\begin{defn} An element $a\in F^k$ which satisfies
$$
\qquad\qquad da = \phi - r \hskip60pt({\bf The\ spark\ equation})
$$
where $\phi\in E^{k+1}$ and $r\in I^{k+1}$ is a \textbf{homological spark}
of degree $k$. Two homological sparks are \textbf{equivalent} if their
difference is of the form
$$
db+s
$$
with $b\in F^{k-1}$ and $s\in I^k$. Given a homological spark $a\in F^k$,
let $\hat a$ denote the equivalence class containing the spark $a$, and let
$\wh\BH^k$ denote the space of homological spark classes.
\end{defn}

Sometimes it is convenient to set $\wh\BH^{-1}=H^0(I^*)$. (The spark
equation is $0=\phi-r$.)

\begin{lem} Each homological spark $a\in F^k$ uniquely determines $\phi\in
E^{k+1}$ and $r\in I^{k+1}$, and $d\phi=0$, $dr=0$. Moreover, $\d_1\hat
a=\phi\in E^{k+1}$ and $\d_2\hat a=[r]\in H^{k+1}(I^*)$ only depend on the
spark class $\hat a\in\wh\BH^k$.
\end{lem}

\begin{proof} Uniqueness of $\phi$ and $r$ is immediate from Axiom (B).
Since $d\phi=dr$, Axiom (B) implies both must vanish. Changing the
homological spark $a$ by adding $db$ does not change $\phi$ or $r$.
Changing the homological spark $a$ by adding $s\in I^k$ does not change
$\phi$ and changes $r$ to $r-ds$ (Again using (B)).
\end{proof}

Given a spark class $c\in\wh\BH^k$,  the element $\phi=\d_1 c\in E^{k+1}$
will be referred to as the {\bf curvature} of $c$, and the class
$\d_2 c\in H^{k+1}(I^*)$ will be referred to as the {\bf divisor class} of
$c$.

Let $Z^k_I(E^*)$ denote the space of cycles $\phi\in E^k$ which are
$F^*$-homologous to some $r\in I^k$.  Let us also define
$$
H^*_I(F^*) \equiv \text{Image}\{H^*(I^*)\to H^*(F^*)\}
\quad\text{and}\quad
H^*(F^*,I^*) \equiv \text{Ker}\{H^*(I^*)\to H^*(F^*)\}.
$$
Finally, let $\wh\BH^k_E$ denote the space of spark classes that can be
represented by a homological spark $a\in E^k$. Note that
$$
\wh\BH^k_E\cong E^k/Z^k_I
$$

\begin{prop} The following diagram commutes, and each row and column is
exact $(k>0)$:
\begin{equation*}
\begin{CD}
\ @.  0 @.  0 @.  0 @.  \ \\
@.    @VVV    @VVV    @VVV @. \\
0  @>>>  \frac{H^k(F^*)}{H^k_{I}
(F^*)}   @>>>  \wh\BH^k_E
@>{\d_1}>> d E^k   @>>>   0  \\
@.    @VVV    @VVV    @VVV @. \\
0  @>>> H^k(F^*/I^*)  @>>>  \wh\BH^k
@>{\d_1}>> Z^{k+1}_I   @>>>   0  \\
@.    @V{\d_2}VV    @V{\d_2}VV    @VVV @. \\
0  @>>>  H^{k+1}(F^*,I^*)  @>>> H^{k+1}(I^*)
@>{}>>  H^{k+1}_{I}(F^*)  @>>>   0  \\
@.    @VVV    @VVV    @VVV @. \\
\ @.  0 @.  0 @.  0 @.  \
\end{CD}
\end{equation*}
\end{prop}

The proof is straightforward and omitted (see [HLZ]).

A spark complex
\begin{align*}
E^*\subset &F^*\\
&\cup\\
&I^*
\end{align*}
is said to be a \textbf{subspark complex} of another spark complex
\begin{align*}
\ov E^*\subset &\ov F^*\\
&\cup\\
&\ov I^*
\end{align*}
if $F^*\subset\ov F^*$ is a differential subcomplex with $E^*=\ov E^*$ the
same,
and $I^*\subset\ov I^*$ such that:
$$
H^*(I^*)\cong H^*(\ov I^*).
$$
Two spark complexes are said to be \textbf{compatible} if they can be
embedded as subspark complexes of (the same) spark complex.

\begin{prop} Given two compatible spark complexes there is a natural
isomorphism $\wh\BH^k\cong\wh{\ov {\bold H}}^k$, $k=-1$, $0$, $1,\dots$.
\end{prop}

The following Lemma will be used twice in the proof of this Proposition.

\begin{lem} Suppose $F^*$, $d$ is a subdifferential complex of $\ov F^*$,
$d$. Then the assertion:
$$
H^k(F^*)\cong H^k(\ov F^*)\qquad k=0,1,\dots
$$
can be restated as follows: Given $g\in F^{p+1}$ and a solution $\a\in\ov
F^p$ to the equation $d\a=g$ there exists $\g\in\ov F^{p-1}$, with
$$
f = \a + d\g\in F^p.
$$
(Not only is there a solution $f\in F^p$ to $df=g$ but one which is
homologous to the given solution $\a$ in $\ov F^p$.)
\end{lem}

\begin{proof} Since $H^{p+1}(F^*)\to H^{p+1}(\ov F^*)$ is injective, there
exists a solution $h\in F^p$ with $dh=g$. Thus $\a-h$ is $d$-closed. Since
$H^p(F^*)\to H^p(\ov F^*)$ is surjective, there exist $\g\in\ov F^{p-1}$ such
that $\a-h+d\g\in F^p$.
\end{proof}

\begin{proof}[\textbf{Proof of Proposition}] We may assume that a subspark
complex is given. Suppose $\bar a\in\ov F^k$ is an $(\ov F, \ov I)$-spark,
i.e. $d\bar a=\phi-\bar r$ with $\phi\in E^{k+1}$ and $\bar r\in\ov
I^{k+1}$. Since $dr=0$ and $H^{k+1}(F^*)\cong H^{k+1}(\ov F^*)$ there exist
$\bar s\in\ov I^k$, $r\in I^{k+1}$ with $\bar r=r-d\bar s$. Therefore,
$d(\bar a-\bar s)=\phi-r$. Since $\phi-r\in F^{k+1}$ and $H^*(F^*)\cong
H^*(\ov F^*)$, the Lemma implies that there exist $\bar b\in\ov F$ with
$a=\bar a-\bar s-d\bar b\in F^k$. This proves that the $(\ov F, \ov
I)$-spark $\bar a\in F^k$ is equivalent to an $(F,I)$-spark $a\in F^k$.

Suppose that $a\in F$ is an $(F,I)$-spark which is equivalent to zero as an
$(\ov F, \ov I)$-spark i.e. $da=\phi-r$, with $\phi\in E^{k+1}$, $r\in
I^{k+1}$ and $a=d\bar b+\bar s$ with $\bar b\in\ov F^{k-1}$ and $\bar
s\in\ov I$. Then $d\bar s=\phi-r$ which implies that both $\phi=0$ and
$d\bar s=-r$. Applying the Lemma to $I^*\subset\ov I^*$ there exist $\bar
t\in\ov I^{k-1}$, $s\in I^k$ with $s=\bar s+d\bar t$. Therefore $a-s=d(\bar
b-\bar t)$. Since $H^k(F^*)\cong H^k(\ov F^*)$, there exist $b\in F^{k-1}$
with $a-s=db$. That is, $a$ is $F$, $I$ equivalent to zero.
\end{proof}

\noindent
{\bf Remark 1.6}  Proposition 1.4 can be strengthened by weakening the
notion of compatible spark complexes to include the case of a chain map from
$F^* \to \overline{F}^*$.  See Section 10 for an application.
\medskip

The remainder of this paper will be devoted to the study of a particular
equivalence class of spark complexes, which will be collectively called
{\bf $(\BR/\BZ)$-spark complexes}. The $3\times 3$ grid for these complexes
is given at the end of the next section.

\section{Sparks (DeRham-Federer)} For completeness we review in this section
some of the results of [HLZ].

\begin{defn} A \textbf{spark of degree $k$} is a current $\a\in\sd k(X)$
with the property that
\begin{equation}
da = \phi - R,\hskip60pt({\bf The\ spark\ equation})\tag*{(2.2)}
\end{equation}
where $\phi\in\SE^{k+1}(X)$ is smooth and $R\in I\SF^{k+1}(X)$ is integrally
flat. Let $\SS^k(X)$ denote the space of sparks of degree $k$.
\end{defn}

Recall that a current $R$ is \textsl{integrally flat} if it can be written as
$R=T+dS$ where $T$ and $S$ are locally rectifiable currents.
The \textbf{DeRham-Federer spark complex} is obtained by taking $E^k=\SE^k(X)$,
$F^k=\sd k(X)$ and $I^k=I\SF^k(X)$. Condition (A) is, of course, standard
[deR].
A proof of condition (B), i.e., that
\begin{equation}
\SE^k(X)\cap I\SF^k(X)=\{0\}\quad k=1, 2,\dots\tag*{(2.3)}
\end{equation}
is given in [HLZ] following Lemma 1.3.

\addtocounter{thm}{2}
\begin{defn}(\textbf{Spark Equivalence}) Two sparks $a$ and $a'$ are
\textbf{(spark) equivalent} if there exists $b\in\sdk(X)$, an arbitrary
current, and $\SS\in I\SF^k(X)$ an integrally flat current, with
\begin{equation}
a - a' = db + S.\tag{2.5}
\end{equation}

The equivalence class determined by a spark $a\in S^k(X)$ will be denoted
$\hat a$, and the space of spark classes will be denoted by
$\wh\BH^k_{\spark}(X)$.
\end{defn}

There are many useful subspark complexes of the DeRham-Federer spark complex.
For example, let $\loc(X)^k$ denote the space of currents of degree $k$ which
can be expressed as differential forms with locally Lebesgue integrable
coefficients, and let $\CDR^k(X)$ denote the space of currents (the
\textbf{current chains} of deRham) that can locally be expressed as
integration over a smooth singular $(n-k)$-chain with integer coefficients.
Take $F^k=\loc(X)^k+d\loc(X)^{k-1}$. Then $I^k=\CDR^k(X)\subset F^k$ and one
obtains $\wh\BH^k_{\text{spark}}(X)\cong\{a\in\loc(X)^k:da=\phi-R$ with
$\phi\in\SE^{k+1}(X)$ and $R\in \CDR^{k+1}(X)$ a current
chain\}/$d\loc(X)^k+\CDR^k(X)$. All such modifications follow from
Proposition 1.4.

Since $H^k(\SE^*(X))=H^k(X, \BR)$, $H^k(I\SF^*(X))=H^k(X, \BZ)$, and
$H^k(\sd *(X)/I\SF^*(X))=H^k(X,S^1)$ (see [HLZ] and [HZ$_1$]), the diagram in
Proposition 1.3 can be rewritten as:
\begin{equation*}
\begin{CD}
\ @.  0 @.  0 @.  0 @.  \ \\
@.    @VVV    @VVV    @VVV @. \\
0  @>>>  \frac{H^k(X,\BR)}{H^k_{\free}(X,\BZ)}   @>>>  \wh\BH^k_\infty(X)
@>{\d_1}>> d\SE^k(X)   @>>>   0  \\
@.    @VVV    @VVV    @VVV @. \\
0  @>>> H^k(X,S^1)  @>>>  \wh\BH^k(X)
@>{\d_1}>> \SZ^{k+1}_0(X)   @>>>   0  \\
@.    @V{\d_2}VV    @V{\d_2}VV    @VVV @. \\
0  @>>>  H^{k+1}_{\tor}(X,\BZ)  @>>> H^{k+1}(X,\BZ)
@>{}>>  H^{k+1}_{\free}(X,\BZ)  @>>>   0  \\
@.    @VVV    @VVV    @VVV @. \\
\ @.  0 @.  0 @.  0 @.  \
\end{CD}
\end{equation*}

\medskip

\section{Hypersparks} Suppose $\SU=\{U_i\}$ is an open covering of $X$ (with
each intersection $U_I$ contractible). Consider the \textbf{\v Cech-Current
bicomplex}
\begin{equation}
\bigoplus_{p,q\ge0}C^p(\SU,\sd q)
\end{equation}
with \textbf{total differential} $D=(-1)^q\d+d$.

Note that $D^2=0$, and that both the horizontal and vertical homology is
zero except along the two edges.

The kernel of $\d$ on the left (vertical) edge $(p=0)$ is the deRham complex
$\sd *(X)$,~$d$ since
$$
0\lra\sd q(X)\lra C^0(\SU,\sd q)\overset{\d}{\lra}C^1(\SU,\sd q)
$$
is exact.

The kernel of $d$ on the bottom edge $(q=0)$ is the \v Cech complex
$C^p(\SU,\BR)$, $\d$ since
$$
0\lra C^p(\SU,\BR)\lra C^p(\SU,\sd o)\overset{d}{\lra}C^p(\SU,\sd 1)
$$
is exact.

We will consider $C^p(\SU,\BZ)$ as a subcomplex of $C^p(\SU,\BR)$, and
$C^p(\SU,IF^q)$ as a subbicomplex of $C^p(\SU,\sd q)$.

\setcounter{thm}1
\begin{defn} A \textbf{hyperspark of degree $k$} is an element
$$
A\in\bigoplus_{p+q=k}C^p(\SU,\sd q)
$$
with the property that
\begin{equation}
D A = \phi - R\hskip60pt({\bf Hyperspark\ equation})\tag*{(3.3) }
\end{equation}
where $\phi\in\SE^{k+1}(X)$ is smooth of bidegree $0$, $k+1$ and
$R\in\bigoplus_{p+q=k+1}C^p(\SU,IF^q)$ is an integrally flat cochain.
\end{defn}

Said differently, we are defining a spark complex by letting
$\overline{F}^*=\bigoplus_{p+q=*}C^p(\SU,\sd q)$ be the \v Cech-current
bi-complex, and setting $\overline{E}^*=\SE^*(X)\subset C^0(\SU, \sd q)$ and
$\overline{I}^*=\bigoplus_{p+q=*}C^p(\SU,IF^q)$. Axioms A and B follow easily.

\setcounter{thm}3
\begin{defn} (\textbf{Equivalence}) Two hypersparks $A$ and $\ov A$ are said to
be \textbf{equivalent} if there exists $B\in\bigoplus_{p+q=k-1}C^p(\SU,\sd
q)$ and
$S\in\bigoplus_{p+q=k}C^p(\SU,IF^q)$ satisfying
\begin{equation}
A - A' = D B + S.\tag*{(3.5)}
\end{equation}
\end{defn}

The equivalence class determined by a hyperspark $A$ will be denoted by
$\wh A$,
and the space of hyperspark classes will be denoted by
$\wh\BH^k_{\hypers}(X)$.

There is obviously a well defined map
$$
\wh\BH^k_{\spark}(X)\lra\wh\BH^k_{\hypers}(X)
$$

\setcounter{thm}5
\begin{thm}
$$
\wh\BH^k_{\spark}(X)\cong\wh\BH^k_{\hypers}(X).
$$
\end{thm}

\begin{proof} The hyperspark complex
$(\overline{F}^*,\overline{E}^*,\overline{I}^*)$  defined above
contains the sparks: $F^*= \sd * (X)$, $E^*=\SE^*(X)$ and $I^*=IF^*(X)$
as a subspark complex (in the component $C^0(\SU, \sd *)$).
Applying Proposition 1.4 completes the proof.
\end{proof}

\section{Smooth Hypersparks} If a spark $a\in S^k(X)$ is smooth then $da=\phi$,
i.e. $R=0$. In fact, $a$ is equivalent to a smooth spark if and only if the
class
of $R$ in $H^{k+1}(X,\BZ)$ is zero. In particular, a general spark need not be
equivalent to a smooth spark. The situation is different for hypersparks.
Consider the \v Cech-form sub-bicomplex of the full \v
Cech-current bicomplex
\begin{equation}
\bigoplus_{p,q\ge0}C^p(\SU,\SE^q)\subset\bigoplus_{p,q\ge0}C^p(\SU,\sd q).
\end{equation}

\addtocounter{thm}{1}
\begin{defn} A hyperspark $A$ is \textbf{smooth} if
$A\in\bigoplus_{p+q=k}C^p(\SU,\SE^q)$.

\bigskip
Note that
\begin{equation}
C^p(\SU,\SE^q)\cap C^p(\SU,IF^q) = \{0\}\tag*{(4.3)}
\end{equation}
unless $q=0$, in which case:
\begin{equation}
C^p(\SU,\SE^0)\cap C^p(\SU,IF^0) = C^p(\SU,\BZ).\tag*{(4.4)}
\end{equation}

Consequently, each smooth hyperspark $A$ satisfies
\begin{equation}
D A = \phi - R\tag*{(4.5)}
\end{equation}
with $R\in C^{k+1}(\SU,\BZ)$ a \v Cech integer cocycle.
\end{defn}

Said differently we are considering the subcomplex with
$F^k=\bigoplus_{p+q=k}C^p(\SU,\SE^q)$,
$E^k=\SE^k(X)\subset C^0(\SU,\SE^k)$, and
$I^k= C^k(\SU,\BZ) \subset C^k(\SU,\SE^0)$.

\begin{align*}
&\SE^2(X)\subset C^0(\SU,\SE^2)\lra \\
&\quad\ \ \uparrow d\hskip25pt\uparrow d   \hskip50pt\uparrow\\
&\SE^1(X)\subset C^0(\SU,\SE^1)\overset{-\d}{\lra} C^1(\SU,\SE^1)\lra\\
&\quad\ \ \uparrow d  \hskip25pt\uparrow d   \hskip50pt\uparrow d  
\hskip45pt\uparrow\\
&\SE^0(X)\subset C^0(\SU,\SE^0)\overset{\d}{\lra} C^1(\SU,\SE^0)
\overset{\d}{\lra} C^2(\SU,\SE^0)\lra\\
&\hskip55pt\cup\hskip60pt\cup\hskip50pt\cup\\
&\hskip45pt C^0(\SU,\BZ)\overset{\d}{\lra} C^1(\SU,\BZ)\overset{\d}{\lra}
C^2(\SU,\BZ)\lra
\end{align*}
\centerline{\bf The Smooth \v Cech-deRham Bicomplex}

\bigskip
\noindent\textbf{Low Degree Cases}:

\medskip
\textbf{Degree 0}: $A\in C^0(\SU,\SE^0)$, $\phi\in\SE^1(X)$, $R\in
C^1(\SU,\BZ)$
\begin{align*}
dA_\a &= \phi\\
A_\b - A_\a &= -R_{\a\b}.
\end{align*}

\textbf{Degree 1}: $A\in C^0(\SU,\SE^1)\op C^1(\SU,\SE^0)$, $\phi\in\SE^2(X)$,
$R\in C^2(\SU,\BZ)$
\begin{align*}
dA_\a &= \phi\\
A_\b - A_\a &= -dA_{\a\b}\\
A_{\b\g} - A_{\a\g} + A_{\a\b} &= -R_{\a\b\g}
\end{align*}

\textbf{Degree 2}: $A\in C^0(\SU,\SE^2)\op C^1(\SU,\SE^1)\op C^2(\SU,\SE^0)$,
$\phi\in\SE^3(X)$, $R\in C^3(\SU,\BZ)$
\begin{align*}
dA_\a &= \phi\in\SE^3(X)\\
A_\b - A_\a &= dA_{\a\b}\\
A_{\a\b} + A_{\b\g} + A_{\g\a} &= dA_{\a\b\g}\\
A_{\b\g\d} - A_{\a\g\d} + A_{\a\b\d} - A_{\a\b\g}  &= -R_{\a\b\g\d}\in\BZ
\end{align*}

\textbf{Degree k}:
$A=\bigoplus_{p+q=k}A^{p,q}\in\bigoplus_{p+q=k}C^p(\SU,\SE^q)$,
$\phi\in\SE^{k+1}(X)$, $R\in C^{k+1}(\SU,\BZ)$
\begin{align*}
dA^{0,k} &= \phi\\
dA^{1,k-1} + (-1)^k\d A^{0,k}  &= 0\\
&\vdots\\
dA^{k,0} - \d A^{k-1,1}   &= 0\\
\d A^{k,0}  &= -R.
\end{align*}

For smooth hypersparks it is natural to define equivalence using only \v
Cech-deRham forms.

\addtocounter{thm}{3}
\begin{defn} Two smooth hypersparks, $A$ and $\ov A$, of degree $k$, are
said to
be \textbf{smoothly equivalent} if there exist
$B\in\bigoplus_{p+q=k-1}C^p(\SU,\SE^q)$ and $S\in C^k(\SU,\BZ)$ satisfying
\begin{equation*}
A - \ov A = D B + S.
\end{equation*}

\bigskip
\noindent\textbf{Low Degree Cases}:

\noindent\textbf{Degree 0}: $A_\a-\ov A_\a=S_\a\in C^0(\SU,\BZ)$.

\medskip
\noindent\textbf{Degree 1}: 
\begin{align*}
A_\a-\ov A_\a &= dB_\a\\
A_{\a\b}-\ov A_{\a\b} = B_\b - B_\a + S_{\a\b}
\end{align*}
for some $B\in C^0(\SU,\SE^0)$ and some $S\in C^1(\SU,\BZ)$.

\noindent\textbf{Degree 2}: 
\begin{align*}
A_\a-\ov A_\a &= dB_\a\\
A_{\a\b}-\ov A_{\a\b} = -B_\b + B_\a + dB_{\a\b}\\
A_{\a\b\g}-\ov A_{\a\b\g} = B_{\b\g} - B_{\a\g} + B_{\a\b} + S_{\a\b\g}
\end{align*}
for some $B\in C^0(\SU,\SE^1)\op C^1(\SU,\SE^0)$ and some $S\in
C^2(\SU,\BZ)$.

Let $\wh\BH^k_{\sm}(X)$ denote the space of smooth hypersparks under smooth
equivalence. There is a natural map
$$
\wh\BH^k_{\sm}(X)\lra\wh\BH^k_{\hypers}(X).
$$
\end{defn}

\begin{thm}
$$
\wh\BH^k_{\sm}(X)\cong\wh\BH^k_{\hypers}(X).
$$
\end{thm}

\begin{proof} Note that the smooth hypersparks form a subspark complex
of the hypersparks and apply Proposition 1.4.
\end{proof}

\section{Grundles} In this section we present grundles which are
closely related (and equivalent) to smooth hypersparks but more geometric in
nature. A grundle of degree $k$ is obtained from a smooth hyperspark $A$ of
degree $k$ by simply replacing the last component $A^{k,0}$ with its
exponential:
\begin{equation}
g_{\alpha_0...\alpha_k}\ =\ e^{2\pi i A_{\alpha_0...\alpha_k}}
\tag{5.1}
\end{equation}
Thus $g\in C^k({\SU}, {\SE}_{S^1})$ where ${\SE}_{S^1}$ is the
sheaf of smooth $S^1$-valued functions. The last component of the
hyperspark equation:
$$
\d A^{k,0} = -R
$$
where $R_{\a_0\dots\a_k}\in\BZ$ implies that
$$
\d g = 0,
$$
i.e. $g$ is a cocycle.

\setcounter{thm}1
\begin{defn} A \textbf{grundle of degree $k$} is a pair $(A,g)$ with $g\in
C^k(\SU,\SE_{S^1})$, $\d g=1$ (a cocycle) and $A\in\bigoplus_{\substack{
p+q=k\\(p,q)\neq(k,0)}}C^p(\SU,E^q)$, satisfying
\begin{align*}
dA^{0,k} &= \phi\in\SE^{k+1}(X)\\
dA^{1,k-1} &+ (-1)^k\d A^{0,k} = 0\\
&\vdots\\
dA^{k-1,1} &+ \d A^{k-2,2} = 0\\
\f1{2\pi i}\f{dg}g &- dA^{k-1,1} = 0
\end{align*}
\end{defn}

The cocycle condition $\d g=1$ implies that for any choice
$A^{0,k}=\f1{2\pi i}\log g$ of the logarithm of the $g$, we have $\d
A^{0,k}\in C^{k+1}(\SU,\BZ)$. Therefore, if we let $R$ denote $-\d
A^{0,k}\in Z^{k+1}(\SU,\BZ)$, then although $R$ is not uniquely determined
by $g$, the class $[R]\in H^{k+1}(\SU,\BZ)$ is  uniquely determined by $g$.
That is, two different choices $A^{0,k}$ and $\ov A^{0,k}$ differ by
$A^{0,k}-\ov A^{0,k}=S\in C^{0,k}(\SU,\BZ)$ and hence $\ov R-R=\d S$
represent the same class in $H^{k+1}(\SU,\BZ)$. In other words
$$
H^k(\SU,\SE_{S^1})\cong H^{k+1}(\SU,\BZ)\qquad(k\ge1)
$$
with $g$ representing the class in $H^k(\SU,\SE_{S^1})$ and $R=-\delta
A^{0,k}$ representing the class in $H^{k+1}(\SU,\BZ)$.  For $k=0$ the
sequence
$$
0\lra\BZ\lra\SE(X)\lra\SE_{S^1}(X)\lra H^1(\SU,\BZ)\lra0\text{ is exact.}
$$

\smallskip
\noindent
\textbf{Low Degree Cases:}

\smallskip
\textbf{Degree 0} A grundle of degree $0$ is a $g\in C^0(\SU,\SE_{S^1})$
satisfying the cocycle condition $\d g=1$. That is $g\in\SE_{S^1}(X)$ is just a
smooth circle valued function.

\bigskip
\textbf{Degree 1} A grundle of degree $1$ consist of a cocycle $g\in
Z^1(\SU,\SE_{S^1})$ (i.e. $g_{ij}\in\SE_{S^1}(\SU_{ij})$ satisfying
$g_{ij}g_{jk}g_{ki}=1$.) and $A=A^{0,1}\in C^0(\SU,\SE^1)$ satisfying
\begin{align*}
&dA_i = \phi\in\SE^2(X)\\
&A_j - A_i  = \f{-1}{2\pi i}\f{dg_{ij}}{g_{ij}}.
\end{align*}
Thus a grundle of degree $1$ is just a hermitian line bundle equipped with a
local trivialization on each $\SU_i$ and equipped with a unitary connection.

\bigskip
\textbf{Degree 2} A grundle of degree $2$ consists of a cocycle $g\in
Z^2(\SU,\SE_{S^1})$ (i.e. $g_{ijk}\in\SE_{S^1}(\SU_{ijk})$ satisfying $\d
g=1$)
$A^{1,1}=\{A_{ij}\}\in C^1(\SU,\SE^1)$, and $A^{0,2}=\{A_i\}\in C^0(\SU,\SE^2)$
satisfying
\begin{align*}
&dA_i = \phi\in\SE^3(X)\\
&A_j - A_i  = dA_{ij}\\
&A_{ij} + A_{jk} + A_{ki} = \f1{2\pi i}\f{dg_{ijk}}{g_{ijk}}.
\end{align*}
Thus a grundle of degree $2$ is just a gerbe with connection (see Hitchin
[Hi]).

\begin{defn} \textbf{Gauge Equivalence} Two grundles $(A,g)$ and $(\ov
A,\bar g)$ are said to be \textbf{gauge equivalent} if: First, there exist an
$h\in C^{k-1}(\SU,\SE_{S^1})$ with $g\bar g^{-1}=\d h$ (i.e. $g$ and $\bar g$
equivalent, via $h$, in $H^k(\SU,\SE_{S^1})$. Second, there exist
$B\in\bigoplus_{\substack{p+q=k-1\\ q\neq0}}C^p(\SU,\SE^q)$ satisfying
the equations
\begin{align*}
A^{0,k} - \ov A^{0,k} &= dB^{0,k-1}\\
A^{1,k-1} - \ov A^{1,k-1}  &= dB^{1,k-2} + (-1)^{k-1}\d B^{0,k-1}\\
&\vdots\\
A^{k-1,1} - \ov A^{k-1,1} &= \f1{2\pi i}\f{dh}h - \d B^{k-2,1},
\end{align*}
\end{defn}

\medskip
\noindent
\textbf{Low Degree Cases}

\smallskip
\textbf{Degree 0} Two grundles $g$ and $\bar g\in\SE_{S^1}(X)$ are gauge
equivalent if they are equal.

\bigskip
\textbf{Degree 1} Two grundles $(A,g)$ and $(\ov A,\bar g)$ are gauge
equivalent
if there exist $h\in C^0(\SU,\SE_{S^1})$ with $g_{ij}{\bar g}^{-1}_{ij}
=h_i/h_j$
and $A_i-\ov A_i=\f1{2\pi i}\f{dh_i}{h_i}$.
Thus equivalence classes of degree $1$ grundles are the same as equivalence
classes of unitary line bundles with unitary connection.

\bigskip
\textbf{Degree 2} Grundles $(A,g)$ and $(A,\bar g)$ (i.e. gerbes with
connection) are gauge equivalent if there exist a cochain $h\in
C^1(\SU,\SE_{S^1})$   and $B\in C^1(\SU,\SE^1)$
satisfying
\begin{align*}
A_i - \ov A_i &= dB_i\\
A_{ij} - \ov A_{ij} &= \f1{2\pi i}\f{dh_{ij}}{h_{ij}} - B_j + B_i\\
g_{ijk}{\bar g}_{ijk}^{-1} &= h_{ij}h_{jk}h_{ki}
\end{align*}

\noindent
{\bf Remark 5.4.}  Let $\wh\BH_{\text{grundle}}^k(X)$ denote the
gauge equivalence classes of degree $k$ grundles on $X$.  It is
straightforward to check that there is an isomorphism
\begin{equation}
\wh\BH_{\text{smooth}}^k(X) \ \cong\ \wh\BH_{\text{grundle}}^k(X)
\tag{5.5}
\end{equation}
induced by exponentiation of the last term as in (5.1).

\section{Cheeger-Simons Sparks} Consider smooth singular (integral) chains
$C_q(X)$ on $X$ and let $C_{\bold Z}^q(X)\linebreak=\text{Hom}(C_q(X),
{\bold Z})$ denote the space of integer cochains of degree $q$. Let
$C^q_\BR(X)\supset C_{\bold Z}^q(X)$ denote the vector space of real
cochains on $X$. Each smooth form $\phi\in\SE^q(X)$ determines, via
integration, a real cochain. This cochain uniquely determines the form
$\phi$, i.e. if a form $\phi$ integrates to zero over all singular chains,
then $\phi=0$. Also, it is easy to see that $\SE^q(X)\cap C_{\bold
Z}^q(X)=0$, $q=1$, $2,\dots$ and
$\SE^0(X)\cap C_{\bold Z}^0(X)=\BZ$. Thus $E^* ={\SE}^*(X)$, $F^*=C_{\bold
R}^*(X)$ and $I^*= C_{\bold Z}^*(X)$ form a spark complex. The associated
sparks are called ``Cheeger-Simons sparks'', that is we have the following.

\begin{defn} A \textbf{Cheeger-Simons spark} of degree $k$ on $X$ is a real
cochain $a\in C_{\bold R}^*(X)$ with the property that
$$
da = \phi - r
$$
with $\phi\in\SE^{k+1}(X)$ a smooth form and $r\in C^{k+1}_\BZ(X)$ an (integer)
cochain. Two such sparks are \textbf{equivalent} if they differ by a real
cochain of the form $db+s$ where $b\in C^{k-1}_\BR(X)$ and $s\in C^k_\BZ(X)$.
Let $\wh\BH^k_{\text{CS}}$ denote the space of Cheeger-Simons spark classes.
\end{defn}

\begin{thm} There exists a natural isomorphism
$$
\wh\BH^k_{\sm}(X)\cong \wh\BH^k_{\text{CS}}(X).
$$
\end{thm}

\begin{proof} We must show that the spark complex of smooth hypersparks
$$
E^*=\SE^*(X),\ F^* = \bigoplus_{p+q=*}C^p(\SU,\SE^q),\text{ and }
I^*=C^*(\SU,\BZ)
$$
and  the spark complex of Cheeger-Simons sparks
$$
E^* = \SE^*(X), F^*=C^*_\BR(X),\text{ and } I^*=C^*_\BZ(X)
$$
are compatible. Let $C^q_\BR$ denote the sheaf of germs of real
$q$-cochains and
let $C^q_\BZ$ denote the subsheaf of (integer) cochains. The spark complex with
\begin{align*}
\ov E^* &=\SE^*(X), \\
\ov F^* &= \bigoplus_{p+q=*}C^p(\SU,\SC^q_\BR),\quad\text{and}\\
\ov I^* &= \bigoplus_{p+q=*}C^p(\SU,\SC^q_\BZ)
\end{align*}
contains both of the previous spark complexes. This larger ``cochain
hyperspark complex'' and Proposition 1.4 complete the proof.
\end{proof}

\section{Differential Characters} In their fundamental paper [CS]   Cheeger
and Simons defined the space of differential characters of degree $k$ to be
the group
$$
\text{DiffChar}^k(X) = \{h\in \hom(Z_k(X), \BR/\BZ) : d h
\equiv\phi\ \ (\text{mod} \ \ \BZ) \text{ for some }\phi\in
\SE^{k+1}(X)\}.
$$
where  $Z_k(X) \subset C_k(X)$ denotes the group of smooth singular
$k$-cycles with $\BZ$-coefficients.

\begin{prop} For any manifold $X$ there is a natural isomorphism
$$
\wh{\BH}_{\text{\rm CS}}^k(X)\ \cong\ \rm{DiffChar}^k(X)
$$
\end{prop}

\begin{proof}
Suppose $a\in C^k_{\BR}(X)$ is a Cheeger-Simons spark.  Define $h_a$ to be
$a$ restricted to $Z_k(X)$ modulo $\BZ$.  If $a=db +s$ represents the zero
spark class, then $h_a(c)=db(c)+s(c)=s(c)\in \BZ$ for every $c\in Z_k(X)$.
Thus $h_a$ depends  only on the spark class $\wh a$.

\medskip
\noindent
{\bf Onto:}  Now given $h\in \hom(Z_k(X), \BR/\BZ)$ we may lift to a
cochain $a\in C^k_{\BR}(X)$.  Moreover, $d h\equiv \phi
\ \ (\text{mod} \ \ \BZ) $ is equivalent to $da-\phi$ being $\BZ$-valued
on $C_{k+1}(X)$, i.e., $da-\phi = -r \in C^{k+1}_{\BZ}(X)$ is an integer
cochain. This proves that  $a$ is a Cheeger-Simons spark.

\medskip
\noindent
{\bf One-to-One:} Suppose $a\in C^k_{\BR}(X)$ is a Cheeger-Simons spark
(i.e., $da=\phi-r$ with $\phi \in \SE^{k+1}(X) \subset C^{k+1}_{\BR}(X) $
and $r \in C^{k+1}_{\BZ}(X)$.  Let $h_a$ denote the induced differential
character.  Suppose $h_a=0$, that is, suppose
$$
a:Z_k(X)\ \to\ \BZ
$$
is integer-valued. Pick an extension $s\in \hom(C_k(X),\BZ)$, i.e., $s\in
C^{k}_{\BZ}(X)$ and $a-s$ vanishes on $Z_k(X)$. Consequently, $a-s=db$ for
some $b\in C^{k-1}_{\BR}(X)$, i.e., the Cheeger-Simons spark $a$ is
equivalent to zero.  
\end{proof}

\section{Holonomy Maps}
Recall now that every $C^{\infty}$ singular integral $k$-chain $c\in C_k(X)$
determines a current $\wt c\in \sd k(X)$  by integration of forms over $c$.
The image of this map, denoted ${\wt C}_k(X)\subset \sd k(X)$, is the deRham
group of {\sl current chains} in dimension $k$. There is a short exact
sequence of chain complexes
\begin{equation*}
0\to N_k(X) \to { C}_k(X)
\overset{\rho}{\longrightarrow} {\wt C}_k(X)\to 0\tag{8.1}
\end{equation*}
and the map $\rho(c)=\wt c$ induces an isomorphism in homology (cf. [deR]).
The elements of $ N_k(X)$ will be called {\sl null chains}. To get some
feeling for this, note that integration over a compact oriented
$k$-dimensional submanifold $\Sigma\subset X$ defines an element
$[\Sigma]\in {\wt C}_k(X)$.  Every smooth triangulation of $\Sigma$
yields a chain $c\in { C}_k(X)$ with $\rho(c) =[\Sigma]$.

Let ${\wt Z}_k(X)$ denote the group of cycles in ${\wt C}_k(X)$.
We define the set of {\bf holonomy maps} of degree  $k$ to be the
group
$$
{{\bold H}{\bold o}{\bold l}}^k(X)\ =\
\{H\in \hom({\wt Z}_k(X), \BR/\BZ) : d H
\equiv  \phi\ \ (\text{mod} \ \ \BZ) \text{ for some }\phi\in
\SE^{k+1}(X)\}.
$$
The form $\phi$ is called the associated {\it curvature form}.
Here $dH(c)=H(\partial c)$ by definition.  Of course, if $H\in
{{\bold H}{\bold o}{\bold l}}^k(X)$, then $h=H\circ \rho\in
\text{DiffChar}^k(X)$.

\setcounter{thm}1
\begin{thm} Every differential character $h\in\hom(Z_k(X), \BR/\BZ)$
descends to a homomorphism $H\in\hom({\wt Z}_k(X), \BR/\BZ)$, i.e., it
can be written in the form $h=H\circ \rho$. This yields a natural
isomorphism
$$
\rm{DiffChar}^k(X) = \ \ {{\bold H}{\bold o}{\bold l}}^k(X)
$$
\end{thm}

\begin{proof}
The pull-back map  ${{\bold H}{\bold o}{\bold l}}^k(X)\lra
\text{DiffChar}^k(X)$ given by
$H\mapsto H\circ\rho$, is obviously injective. Therefore, the
isomorphism follows from the first assertion. To prove this let $h\in
\text{DiffChar}^k(X)$ be any differential character. By
definition $d h = \phi$ for some smooth $(k+1)$-form $\phi$. Now it
suffices to show that for any pair $c_1, c_2\in Z_k(X)$ satisfying ${\wt
c}_1={\wt c}_2$, there exists a null chain $B$ with  $
c_1-c_2\ =\ dB
$
because in this case
$$
h(c_1)-h(c_2)\ =\ h(dB)\ =\ \int_B \phi \ =\ 0.
$$
Note that $a=c_1-c_2$ is a null cycle. Furthermore, the homology class of $a$
in $C_*(X)$ must be  zero since $\wt a=0$ and the map $\rho(c)\equiv \wt c$
in (7.1) is injective in homology. Hence, there exists $b\in C_{k+1}(X)$
with $db=a$, and we have $d\wt b=\wt a=0$. Now since the map $\rho$ is
surjective in homology, we can write $\wt b = \wt z +d \wt e$ where
$dz=0$. Set $B=b-z-de$. Then $B$  is a null chain with $dB=a=c_1-c_2$ as
desired.
\end{proof}

\section{Grundle holonomy} Combining Remark 5.4, Theorem 8.2, Proposition
7.1 and Theorem 7.3 we see that every grundle $(A,g)$ of degree $k$ induces
a homomorphism
$$
{h}_{(A,g)} : {\wt Z}_k(X) \ \lra \  \BR/\BZ = S^1
$$
called the {\bf holonomy} of $(A,g)$.  It depends only on the
gauge-equivalence class of the grundle.  (By Theorem 4.7 every hyperspark
also has such a holonomy homomorphism.)

\def\Sig{M}
In particular  let $\Sig$ be a  compact connected oriented manifold of
dimension $k$ ($\partial \Sig = \emptyset$) and $f:\Sig \to X$ a smooth map.
Then by Theorem 8.2 the holonomy of $(A,g)$ on $\Sig$
$$
{h}_{(A,g)}(\Sig) \in S^1
$$
is well-defined (i.e., independent of the choice of triangulation
of $\Sig$).  It  can be computed as follows.  The induced grundle $f^*(A,g)$
represents a class in ${\widehat{\bold H}}_{\text{grundle}}^k(\Sig)
=H^k(\Sig;\,\BR)/H^k(\Sig;\,\BZ) = \hom(H_k(\Sig;\,\BZ), S^1)
= \hom(\BZ \cdot [\Sig],   S^1)$
by the diagram at the end of \S 2.  Then ${h}_{(A,g)}(\Sig)=
[f^*(A,g)]([\Sig])$.

What is most interesting is to find intrinsic geometric definitions of the
holonomy over  such cycles $\Sig$ and to find combinatorial formulas for
${h}_{(A,g)}(\Sig)$ in terms of the grundle data. A basic example of this
occurs in degree 1 where equivalence classes of grundles coincide with gauge
equivalence classes of complex line bundles with unitary connection.  Here
the holonomy around a closed loop $\gamma$ coincides with the rotation
angle obtained by parallel translation around $\gamma$ (the classical
holonomy) [CS]. This is expressed in the combinatorial grundle formula
derived below.

We shall systematically derive formulas for grundles of all degrees.

The general procedure is as follows. Consider a grundle of degree
$k$ presented by a smooth hyperspark $A$. (That is, we choose logarithms
$A_{\a_0\dots \a_k}$ for the $g_{\a_0\dots \a_k}$ as in (5.1).)
Let $h$ be a Cheeger-Simons spark which is equivalent to $A$ in the
cochain hyperspark complex.  Then the holonomy ${h}_A:Z_k(X)\to
\BR/\BZ$ is given by
$$
{h}_A(z)\ \equiv\ h(z)\ \ \ \  (\text{mod }  \BZ )
$$
The idea now is to use the equivalence between $A$ and $h$ to generate
the formula.

\medskip
\noindent{\bf Degree 0.}  Here a 0-grundle is given by  $g_\a=\exp(2\pi
iA_\a)$ on $U_\a\in {\SU}$.  Let $h\in C^0_{\BR}(X)$ be an equivalent
Cheeger-Simons spark.  Then $A_\a-h = S_\a$ in $U_\a$ for some $S\in
C^0(\SU, C^0_{\BZ})$.  We conclude that the holonomy on the zero-cycle
$x\in X$ can be written as
$$
{h}_A(x)\ \equiv A_\a(x)\ \ \ \  (\text{mod }  \BZ )
$$
where $x\in U_\a$. Of course, as noted in \S 5,   a grundle of  degree 0 is
just a  smooth map $g:X\to S^1$, and its holonomy at $x$ is seen to be
$g(x)\in S^1$. This is the ``intrinsic'' definition of holonomy in degree 0.

\medskip
\noindent{\bf Degree 1.} Let  $A=(\{A_{\a\b}\}, \{A_\a\})\in
C^1(\SU, \SE^0)\oplus C^0(\SU, \SE^1)$ be a smooth hyperspark representing
a 1-grundle as above, and let $h\in  C^1_{\BR}(X)$ be an equivalent
Cheeger-Simons spark. The elements $A$ and $h$  both lie in
$C^1(\SU, C^0_{\BR})\oplus C^0(\SU, C^1_{\BR})$ where $h=(\{h_{\a,\b}\},
\{h_\a\}) = (0, h|_{U_\a})$. The equivalence of $A$ and $h$ means  there
exist elements $a\in C^0(\SU, C_{\BR}^0)$ and
$S\in C^1(\SU, C^0_{\BZ})\oplus C^0(\SU, C^1_{\BZ})$ with
\begin{equation}
\begin{split}
da_\a\ &=\ A_\a -h +S_\a\\
a_\b-a_\a\ &=\ \ \  A_{\a\b}  +S_{\a\b}
\end{split}
\end{equation}
Let $\gamma\in \SZ_1(X)$ be a closed loop and write $\gamma = \sum_{\a=1}^N
\gamma_\a$ where $\g_1,...,\g_N$ are successive 1-simplices subordinate
to elements $U_1,...,U_N$ of the covering. Set
$V_{\a,\a+1}=\g_{\a+1}\cap\g_\a$ with indices taken mod $N$.  Then
using Theorem 8.2 and calculating with equation (9.1) we find
\begin{align*}
{h}_A(\g)\ \equiv\ \sum_{\a=1}^N h(\g_\a)\ &\equiv\ \sum_{\a=1}^N
\left\{\int_{\g_\a} A_\a + a_\a(d\g_\a)\right\}\ \ \ \  (\text{mod }  \BZ )\\
&\equiv\ \sum_{\a=1}^N  \left\{ \int_{\g_\a}
A_\a + a_\a(V_{\a,\a+1}-V_{\a-1,\a})\right\}\ \ \ \  (\text{mod }  \BZ )\\
&\equiv\ \sum_{\a=1}^N   \int_{\g_\a}
A_\a -  \sum_{\a=1}^N A_{\a,\a+1}(V_{\a,\a+1})\ \ \ \  (\text{mod }  \BZ ).
\end{align*}

\medskip
\noindent
{\bf Remark 9.2.}  Suppose that our grundle data $\{A_j\}$ and
$\{g_{i,j}\}$, with $g_{ij} =\exp (2\pi i A_{ij})$, correspond to the line
bundle with connection $(L,D)$. Then $A_j$ represents the connection 1-form
with respect to a  nowhere vanishing section $\s_j$ of  $L\bigl|_{U_j}$,
and $g_{ij}$ is the change of trivialization on $U_i\cap U_j$.
One has that $D\s_j= iA_j\otimes \s_j$, and a section $\s = f_j \s_j$ is
{\bf parallel} along the arc $I_j\subset U_j$ iff $df_j+if_j A_j=0$  on
$I_j$.  Therefore the effect of parallel translation in $L$ once around the
loop $\g$, is given by the formula:
\begin{equation*}
{h}_A(\g)\ \equiv\ \sum_{j=1}^N\,\int_{\g_j} A_j \ -\ \sum_{j=1}^N
\int_{V_{j}}A_{j,j+1}\ \ \  (\text{mod }\BZ)\tag{9.3}
\end{equation*}
which we derived above.  The first integrals represent parallel translation
along $I_j$ in the given frame; the second ``integrals'' represent the
change of frame at the vertex.

\medskip
\noindent{\bf Degree 2.} Let  $A=(\{A_{\a\b\g}\},\{A_{\a\b}\}, \{A_\a\})\in
C^2(\SU, \SE^0)\oplus C^1(\SU, \SE^1)\oplus  C^0(\SU, \SE^2)$ be a smooth
hyperspark representing the 2-grundle (or gerbe with connection)  as above,
and let $h\in  C^2_{\BR}(X)$ be an equivalent Cheeger-Simons spark.
Equivalence in the cochain hyperspark complex means that there exist
elements
$a\in C^1(\SU, C_{\BR}^0)\oplus C^0(\SU, C_{\BR}^1)$ and
$S\in C^2(\SU, C^0_{\BZ})\oplus
C^1(\SU, C^1_{\BZ})\oplus C^0(\SU, C^2_{\BZ})$ with
\begin{align*}
da_\a\ \qquad  &=\ A_\a -h +S_\a\\
a_\b-a_\a\ - da_{\a\b}\  &=\ \ \  A_{\a\b}  +S_{\a\b}\tag{9.4}\\
a_{\a\b}+a_{\b\g}+a_{\g\a}\ \  &=\ \ \ A_{\a\b\g} +S_{\a\b\g}
\end{align*}
Suppose now that $\Sigma$ is a compact oriented surface without boundary in
$X$. We assume that  $\Sigma$ is provided with a smooth cell structure $T^*$
dual to some triangulation $T$ of $\Sigma$. The 2-cells of $T^*$ are
polygons defined by taking  the stars of the vertices of $T$ in its  first
barycentric subdivision. Each vertex of $T^*$ meets exactly three edges and
three 2-cells. By taking $T$ fine enough we may assume  each  2-cell $P$ of
$T^*$ is contained in some $U_\a\in\SU$. We choose one such $U_\a \supset P$
and label $P$ as $P_\a$. Each $P_\a$ is oriented by the orientation of
$\Sigma$. Each edge $E$ is contained in exactly two faces $P_\a$ and $P_\b$.
We label $E$ as $E_{\a\b}$ and orient $E_{\a\b}$ as the boundary of
$P_\a$. (Thus $E_{\b\a}$ is oriented oppositely, as the boundary of
$P_\b$.) Each vertex $V$ meets exactly three edges and three faces and is
thereby labelled $V_{\a\b\g}$. Using the equations (9.4) and arguing
as in the degree -1 case, we find that the holonomy of our gerbe with
connection $A$ on $\Sigma$ is given by the formula:
\begin{equation*}
{h}_A(\Sigma)\ \equiv\ \sum_{\a}\,\int_{P_\a} A_{\a} \ -\
\sum_{\a\b}\,\int_{E_{\a\b}} A_{\a\b} \ +\
\sum_{\a\b\g}\,\int_{V_{\a\b\g}}A_{\a\b\g}\ \ \  (\text{mod } \BZ )\tag{9.5}
\end{equation*}

\medskip
\noindent{\bf Degree k.} Let $A\in C^k(\SU, \SE^0)\oplus\dots\oplus
C^0(\SU, \SE^k)$ be a smooth hyperspark and $h\in C^k_{\BR}(X)$ an
equivalent Cheeger-Simons spark. Then as above there exist elements
$a\in C^*(\SU,C^*_{\BR})$ and $S\in C^*(\SU,C^*_{\BZ})$ such that
\begin{equation*}
Da \ =\ A-h+S\tag{9.6}
\end{equation*}
where $D$ is the total differential as in \S 3.  Now let $M$ be a compact
oriented submanifold of dimension $k$ in $X$.  We suppose $M$ is given a
smooth cell structure which is dual to a triangulation and for which each
$k$-cell $e_\a$ is contained in an open set $U_\a\in \SU$. Each
$(k-\ell)$-cell is contained in exactly $\ell+1$ $k$-cells and is thereby
labelled $e_{\a_0\dots\a_\ell}$.  The cell $e_\a$ is oriented by $M$.
Inductively the cell $e_{\a_0\dots\a_\ell}$ is oriented as
part of the boundary of $e_{\a_0\dots\a_{\ell-1}}$.  Then proceeding as
above one finds that
\begin{equation*}
{h}_A(M)\ \equiv\ \sum_{\ell=0}^k(-1)^{\ell}
\sum_{\a_0\dots \a_{\ell}} \int_{e_{\a_0\dots\a_{\ell}}}
A_{\a_0\dots\a_{\ell}} \ \ \  (\text{mod }  \BZ )\tag{9.7}
\end{equation*}
where the $\ell^\supth$ sum is over the $(k-\ell)$-cells in the complex.

\medskip
\noindent
{\bf Remark 9.8. Holonomy via de Rham - Federer Sparks.} Let $A$ and $M$
be as above and suppose $A$ is equivalent to a de Rham - Federer spark
$a\in\SD{'{^k}}(X)$.  By changing $a$ by a boundary $db$, $b\in\sdk(X)$  we
may assume that $a$ is smooth on $M$ (see [HLZ, Prop. 4.2].)  If $a$ is a
dRham - Federer spark which is smooth, then $a$ is also a smooth hyperspark
with only one non-zero component, namely of bidegree $0,k$.  Consequently
(9.7) reduces to:
\begin{equation*}
{h}_A(M)\ \equiv\ \int_M a \ \ \  (\text{mod }  \BZ )\tag{9.9}
\end{equation*}
This direct connection between de Rham - Federer sparks and differential
characters given by (9.9) was the basis for showing
$
\wh{\BH}_{\text{spark}}^k(X)\ \cong\ \wh{\BH}_{\text{Diff Char}}^k(X)
$
in [HLZ].

\section{Further Spark Complexes} There are many interesting homological
spark complexes which are not treated in this paper. We sketch two quick
examples.

\medskip
\noindent
{\bf Combinatorial Sparks.} \
Suppose $X$ is a manifold provided with a smooth triangulation (or
cubulation) $\ST$.  Consider the family of all barycentric subdivisions (or
dyadic subdivisions) of $\ST$, and let $C_k(\ST)$ denote the abelian group
generated by the oriented  $k$-simplices (or $k$-cubes)
belonging to these subdivisions,  modulo the relation
$\sigma =- \sigma^*$ where $^*$ indicates the reversed orientation.
Then we can define
\begin{equation}
F^*=\hom(C_*(\ST), \BR), \qquad E^*=\text{Image}(\SE^*(X)\to F^*),\qquad
I^*=\hom(C_*(\ST), \BZ).
\end{equation}
Note that because of the subdivisions, the differential forms $\SE^*(X)$
inject into $F^*$.  From this and the fact that $C_*(\ST)$ computes the
integral homology of $X$, we have that  $H^*(X;\,\BR)\cong H^*(F^*)\cong
H^*(E^*)$.  Using the subdivisions we also see that $E^k\cap I^k=\{0\}$ for
$k>0$.  Hence, (10.1) is a homological spark complex. Since $H^*(I^*)\cong
H^*(X;\,\BZ)$, one sees that this complex is compatible with the
Cheeger-Simons sparks.  Hence the associated group of homological spark
classes  $\wh{\BH}^k$ coincide with the differential characters
$\wh{\BH}^k_{\text{Diff Char}}(X).$ When $\ST$ is a cubulation, this
gives a direct approach to Cheeger's $*$-product [C].

\medskip
\noindent
{\bf Current Cochain Sparks.} \
Suppose $X$ is a manifold and let ${\wt C}_*(X)$ denote de Rham's complex of
current chains discussed in \S 8. Then we can define the complex of {\sl
current
cochains} by $F^k =  {\wt C}^k(X) = \hom({\wt C}_k(X), \BR)$ with
$ I^k=\hom({\wt C}_k(X), \BZ)$ the {\sl integral current cochains}.
The differential forms $E^*=\SE^*(X)$ clearly inject into $F^*$, and since
${\wt C}_*(X)$ computes the integral homology of $X$, we have that
$H^*(X;\,\BR)\cong H^*(F^*)\cong H^*(E^*)$. Evidently, $E^k\cap I^k=\{0\}$
for $k>0$, and we have  a homological spark complex.

As we saw in section 8, the natural map from $F^*$, current cochains, to
$\overline {F}^*$, cochains, is not injective.  However, as noted in Remark
1.6, Proposition 1.4 can be easily modified to include this case, proving that
current cochain spark class are isomorphic to Cheeger-Simons cochain spark
classes.

Just as it was straightforward to show (cf. Prop. 7.1) that Cheeger-Simons
spark  classes are isomorphic to differential characters, it is
straightforward to show that current cochain spark classes are isomorphic to
${\bold H}{\bold o}{\bold l}(X)$. This provides another proof of Theorem
8.2.

\section{Applications and Examples.}

To give the reader a feeling for the   abundance  and usefulness of spark
complexes we briefly sketch a number of  examples and applications. The
emphasis is on evaluating which of the approaches to the $(\BR,\BZ)$ theory
is best suited for the particular application or example being considered.
Frequently, multiple points of view are enlightening.
\medskip

\noindent
{\bf Ring structure.} One of the deep and most useful features of
$(\BR,\BZ)$-spark classes is the existence of a natural graded  ring
structure\ \   $ * :\wh\BH^*(X)\times \wh\BH^*(X)\to \wh\BH^*(X)$.
The difficulty of establishing  this product depends greatly
 on the formulation (i.e., spark complex) chosen  to represent
characters.  The original
construction  by Cheeger and Simons [CS], [C]
was   quite non-trivial due  to the difference between the
wedge-product of differential forms and the cup-product of the cochains
they define. If one takes a grundle approach or the approach of (even
smooth) hypersparks, the existence of this product seems something of a
miracle.  In the de Rham-Federer spark approach   there is a simple and
quite useful formula for this product ([HLZ], [GS]) which holds for generic
sparks in any pair of spark classes. It is established using transversality
theory for flat and rectifiable currents.

These remarks apply to the  $(\SO,\BZ)$-spark classes discussed in [HL$_5$].

We see that for this ring structure the existence of many different
approaches is quite useful. Another example is the following.

\medskip
\noindent
{\bf Functoriality.}   Any smooth mapping $f:X\to Y$ between manifolds
induces a ring homomorphism $f^*: \wh\BH^*(Y)\to \wh\BH^*(X)$  compatible
with  $\d_1$ and  $\d_2$.  This assertion is evident from say the
Cheeger-Simons or the holonomy approach to characters, and also
in the setting of smooth hypersparks.  However, it is far from clear in
many   other approaches since, for example, the pull-back of currents is not
universally defined.  Here the ability to switch from one theory to the
other is quite useful.

\medskip
\noindent
{\bf Gysin maps and Thom homomorphisms.} To a  smooth proper submersion
$f:X\to Y$ between oriented manifolds there is associated a Gysin mapping
$f_*: \wh\BH^*(X)\to \wh\BH^*(Y)$ compatible with  $\d_1$ and  $\d_2$. This
is completely evident from the de Rham-Federer  viewpoint (cf. [HLZ, \S 10])
but rather mysterious from many other points of view. Similarly using sparks
one can naturally define Thom homomorphisms [HLZ,\S 9] which are not so
evident in other formulations.

\medskip
\noindent
{\bf Secondary (refined) characteristic classes.} Working from their
viewpoint Cheeger and Simons developed a full theory of characteristic
classes associated to princpal bundles with connection [CS]. It
simultaneously refined both Chern-Weil theory and the theory of integral
characteristic classes, and it gave new invariants for flat bundles and
foliations. Their elegant arguments involved functorialty and the existence
of classifying spaces for bundles with connection.

These refined classes can also be defined from the spark viewpoint.
The construction involves a choice of section or bundle map $\a$ and yields
an $L^1_{\text{loc}}$-transgression form $T$ with $dT=\varphi-S(\a)$ where
$\varphi$ is a characteristic form and $S(\a)$ is a rectifiable current
measuring certain singularities of $\a$. The secondary Chern and
Pontrjagiin classes can be systematically developed in this way. The
advantage here is that singularities of geometric problems explicitly enter
the picture. 

\medskip
\noindent
{\bf Example: The Euler Spark.} An understanding of this case provides
geometric motivation for the spark equation and spark equivalence. Suppose
$E$ is a real oriented vector bundle. Equipping $E$ with a metric connection
enables one to compute an Euler (Pfaff) form $\chi$ for $E$. Equipping $E$ with
a section $s$ whose zero set is reasonable (say $s$ has simple zeros or more
generally $s$ is atomic) provides a divisor or zero current Div$(s)$. In
fact, there is a canonically defined $L^1_{\text{loc}}$ spark $\s$
satisfying the local Gauss-Bonnet equation (the spark equation):
$$
d\s = \chi - \Div(s)\qquad\text{on }\ X.
$$
See [HZ$_2$] for the details of this discussion. The explicit formula for $\s$
as a fiber integral leads one to compare two sparks $\s$ and $\s'$ generated
by two sections $s$ and $s'$ by considering the section $\g=ts+(1-t)s'$ of
the bundle $\wt E$, equal to $E$ pulled back over $\BR\times X$. Now a Stokes
theorem argument yields the {\it spark comparison formula}
$$
\s - \s' = dL + R\qquad\text{on }\ X.
$$
where $L$ is a fiber integral and hence an $L^1_{\text{loc}}$ form, and
where $R$ is the pushforward of the divisor of $\g$ (cut off on $0\le
t\le1$), and hence is a rectifiable current.

Note that taking the exterior derivative of both sides yields the {\it
divisor comparison formula}
$$
\Div(s') - \Div(s) = dR.
$$

In summary the refined euler class $\hat{\chi}\in\wh H^n(X)$ of a bundle with
metric connection is the spark class made up of the euler sparks $\s$. It's
divisor class is the euler class $\chi$ of $E$ in $H^n(X,\BZ)$ made up of the
divisors $\Div(s)$, it's curvature is the euler Pfaff form $\chi$, and the
spark equation is just the local Gauss-Bonnet equation above.

\medskip
\noindent
{\bf Flat sparks.} A flat spark is any representative of a class
$c\in\wh\BH^k(X)$ whose curvature form vanishes, i.e., $\phi \in
\ker\,\d_1=H^k(X,S^1)$. The {\sl grundle point of view} is very natural in
this case.  Moreover, the \v Cech representation of $H^k(X,S^1)$ shows that
 there exists a grundle representative of the form $(A,g)$ with $A=0$ and
with $g\in \SZ^k(\SU, S^1)$ actually taking locally constant values in $S^1$.

The {\sl holonomy point of view} is also natural in this case. A holonomy
map is flat if it vanishes on boundaries.  This induces a homomorphism from
$H_k(X;\BZ)\cong H^{n-k}(X;\BZ)$ to $S^1$, so that
${\bold H}{\bold o}{\bold l}^k_{\text{flat}}(X) \cong
\hom(H^{n-k}(X;\BZ),S^1)$. Thus the isomorphism
$$
{\bold H}{\bold o}{\bold l}^k_{\text{flat}} (X)\ \cong\ H^k(X;S^1)
$$
is equivalent to
$$
\hom(H^{n-k}(X;\BZ),S^1)\ \cong\ H^k(X;S^1)
$$
which is one way to state the classical Poincar\'e Duality Theorem
over $\BZ$.  The same discussion applies to  {\sl flat differential
characters}.

The {\sl de Rham-Federer point of view} is natural once the isomorphism
$$
H^k(X;S^1)\ \cong\ \frac{\{a\in\sdk(X) : da\in Rect^{k+1}(X)\}}
{d\sdk(X)+Rect^k(X)}
$$
is established.  (See [HLZ, pg 8] for a proof.) One considers all sparks
$a$ satisfying the spark equation $da=\phi-R$ with $\phi=0$. Note that
spark equivalence is the same as the equivalence relation for $H^1(X,S^1)$.

\medskip
\noindent
{\bf Projective bundles.}  Flat grundles  appear in addressing
the following question. Let $P\to X$ be a smooth fibre bundle with fibre
$\BP^{n-1}$ and structure group PGL$_n(\BC)$. {\sl When does there exist a
complex vector bundle $E\to X$ of rank $n$ with $\BP(E)\cong P$?}
This is equivalent to asking for a complex line bundle $\lambda\to P$ whose
restriction to the fibres is ${\SO}(-1)$, for then $E$ is obtained from
$\lambda$ by blowing down the zero-section.

To answer the question we first reduce the structure group, by general
principles, to PU$_n$. Consider the short exact sequence
$$
0\to S^1 \to U_n\overset{\pi}{\lra} PU_n \to 1.
$$
Let ${\SU}=\{U_i\}_i$ be an acyclic covering of $X$ with trivializations
$P\bigr|_{U_i} \to U_i\times \BP^{n-1}$ and transition functions
$G_{ij}:U_{ij}\to \text{PU}_n$. Choose liftings
$g_{ij}:U_{ij}\to \text{U}_n$ with $\pi\circ g_{ij}=G_{ij}$  and define
a 2-grundle over ${\SU}$ by setting $\phi=0$,  $A_i=0$,   $A_{ij} = \frac
1n d \log \det g_{ij}$ and $g_{ilk}I = g_{ij}g_{jk}g_{ki}$.
Note that $\pi(g_{ij}g_{jk}g_{ki})\equiv1$ and so $g_{ijk}\in S^1$
and that $(A_i,A_{ij},g)$ satisfy the grundle conditions of \S 5.
The class of the cocycle $g_{ijk}$ in $H^2_{\tor}(X;{\SE}_{S^1})\cong
H^3_{\text{tor}}(X;\BZ)$ represents the obstruction to the existence of
the desired vector bundle $E$. This group $H^3_{\text{tor}}(X;\BZ)$,
called the {\sl topological Brauer group}, was introduced by Grothendieck
and Serre [G].

\medskip
\noindent
{\bf Hodge sparks.} On a riemannian manifold the Greens operator
$G$ provides an important source of sparks (see [HLZ, \S 12] for more
details). Given an integrally flat $(k+1)$ cycle $R$ (such as a current
chain without boundary), the {\bf Hodge spark} of $R$ is defined by
$\s(R) =-d^*G(R)$.  It satisfies the spark equation
$$
d\s(R)\ =\ H(R)-R
$$
where $H(R)$ is the harmonic form determined by the class of $R$.

This is an example where the de Rham - Federer point of view is clearly
superior to the other points of view.
The holonomy map is easily computed for any cycle $S$ which does not meet
$R$ by the formula
$$
h_R(S)\ \equiv\ \int_S\,\s(R) \quad  (\text{mod } \BZ)
$$

If $R$ is a boundary, then $H(R)=0$ and the spark class of $\s(R)$ lies in
the torus $H^k(X;\BR)/H^k(X;\BZ)=\ker \d_1 \cap \ker \d_2$.  This
``Jacobian'' torus can be realized either as $Jac^k=Har^k(X)/Har^k_0(X)$,
harmonic forms modulo those with integral periods, or by duality as
$\text{Hom}(Har^k_0(X), \BR/\BZ)$.  The induced ``Abel-Jacobi''map
$j:B^{k+1}\to Jac^k$ on the space of boundaries can be realized directly by
sending the boundary $R=d\Gamma$ to the spark class of $H(\Gamma)$, the
harmonic form corresponding to $\Gamma$.  That is $\s(R)$ and $H(R)$ are
equivalent sparks.

We define boundaries in the kernel of the Abel-Jacobi map $j$ to be {\sl
linearly equivalent to zero}.  These are the ``principal boundaries'' whose
associated spark class is 0. One has that a boundary $R=d\Gamma$ is
linearly equivalent to zero iff
$$
\s(R)\ =\ \int_{\Gamma} \theta \quad  (\text{mod } \BZ)
$$
for all harmonic $(n-k)$-forms $\theta$ with integral periods.

In summary this example is best described from the de Rham-Federer point of
view. See [Hi], [Cha] for a gerbe discussion in special cases.

\medskip
\noindent
{\bf Morse sparks.}  Suppose $\varphi_t$ is a Morse-Stokes flow [HL$_4$]
or more generally a finite-volume flow (see [HL$_4$] and [M]). Then for each
critical point $p$ both the stable manifold $S_p$ and the unstable manifold
$U_p$ have finite volume and define currents (by integration) also denoted
by $S_p$ and $U_p$.  Under the flow $\varphi_t$ each smooth form $\a$ has
the limit
$$
P(\a)\ =\ \lim_{t\to\infty}\varphi^*_t(\a) \ =\ \sum_{p\in \text{Cr}}
\left(\int_{U_p} \,\a\right)\,S_p
$$
in the space of currents.  For a form $\a$ whose residues $\int_{U_p}\a$
are integers, the limit $P(\a)=R$ is a rectifiable current, and $\a$ can be
considered a ``Thom form'' for $\a$.

There exists a continuous linear operator $T$ from forms to currents which
lowers degree by one. The {\bf Morse spark} $T(\a)$  satisfies the
spark equation
$$
dT(\a)\ =\ \a-R.
$$
The operator $T$ is induced by the kernel current $\bold
T\equiv\{(x,\varphi_t(x))\in X\times X:t\geq 0\}$.
There is an abundance of interesting sparks coming from specific flows and
forms.  This is a case where the de Rham - Federer viewpoint seems the only
way to construct the  spark classes.

\medskip
\noindent
{\bf Curvature-driven sparks.}  When $H^k(X,S^1)=0$, the curvature form
uniquely determines each spark class  since
$\d_1: \wh\BH^k(X)\to \SZ_0^{k+1}(X)$ is an isomorphism.  By Poincar\'e
duality
$H^k(X,S^1)$ vanishes if and only if $H_k(X,\BZ)$ vanishes. (Also note
that $H^k(X,S^1)=0$ if and only if $H^k(X,\BR)=0$ and
$H^{k+1}_{\text{tor}}(X,\BZ)=0$.)  In this case the holonomy point of view
is transparent.  Since each cycle $\Sigma$ of dimension $k$ is a boundary
$\Sigma = d\Gamma$, the holonomy $h_{\phi}(\Sigma) \equiv
\int_{\Gamma} \phi\qquad $ (mod $\BZ$) is well defined.  (
$\int_{\Gamma} \phi$ is well defined modulo the periods of $\phi$.)

\medskip
\noindent
{\bf Example: The Wess-Zumino term.}
Suppose $G$ is a compact simple simply-connected Lie group.
Let $\Phi(X,Y,Z) =B(X,[Y,Z])$ denote the Cartan 3-form (bi-invariant and
d-closed) on $G$ where $B$ denotes the Killing form on the Lie algebra
$\ger$.  Since $H^3(G;\BZ)\cong \BZ$, $\Phi$ may be normalized to
be a generator.  Also, $H^2(G;\BZ)=0$, and hence $\Phi$ uniquely determines
a class $wz\in\wh\BH^2(G)$.

As noted above the holonomy map (as well as the differential character)
point of view is transparent in this case.  However, a deRham-Federer spark
also provides some geometric insight. Let $E:\ger\to G$ denote the
exponential map. There exists a unique bounded star-shaped neighborhood $U$
of the origin in $\ger$ which is diffeomorphic, under $E$, to the open set
$G-C$ where $C$ is the {\sl cut locus}.  The cut locus is a codimension-3
stratified set which naturally determines a current of degree 3 on $G$. If
$a\in\SE^2(\ger)$ is chosen to satisfy the equation $da=E^*(\Phi)$ on
$\ger$, then $A=E_*(\chi_U a)$ is an $L^1_{\text{loc}}$-spark on $G$ with
curvature form $\Phi$. In fact, one can show that on $G$
\begin{equation}
dA=\Phi-[C].
\tag{11.1}
\end{equation}
Note that if $\bar A$ is another current on $G$ satisfying the equation
(11.1), then $\bar A= A+dB$ for some $B\in {\SD'}^1(G)$ because
$H^1(G)=0$. In particular $a$, and therefore also $A$, can be chosen to be
Ad$_G$-invariant. Ad$_G$-invariant sparks satisfying (11.1) will be called {\sl
Wess-Zumino sparks}.

An even more explicit construction of an (Ad$_G$-invariant) Wess-Zumino
spark is given by the taking the Hodge spark $A= -d^*G([C])$ of the cut locus
$C$.

\medskip
\noindent
{\bf Torsion classes.}  Suppose $c\in\wh\BH^k(X)$ is a torsion class of
torsion degree $n$.  Then since $nc=0$, the curvature satisfies $n\d_1(c)
= n\phi=0$.  Hence, $c\in\ker\d_1\equiv\wh\BH^k_{\text{flat}}(X)
= H^k(X,S^1)$ is a flat class.  Such classes  $c\in H^k(X,S^1)$ acn be
related to elements in $H^k(X,\BZ_n)$ as follows.  The exact triple
$0 \to \BZ_n\overset{1/n}{\lra} S^1  \overset{n}{\lra} S^1 \to 0$ induces a
long exact sequence $\dots \to  H^k(X,\BZ_n)\overset{1/n}{\lra}
H^k(X,S^1)\overset{n}{\lra} H^k(X,S^1)\lra\dots$. Since $nc=0$, there
exists an element
$u\in H^k(X,\BZ_n)$ whose image in $H^k(X,S^1)$ is $c$.
Let $\b$ denote  the  Bockstein map  induced by
$0\to \BZ\overset{n}{\lra}\BZ \to \BZ_n\to 0$.  Then the ``divisor class''
of $c$, namely $\d_2(c)\in H^{k+1}(X, \BZ)$, can be computed more directly
in terms of $u$.  That is, $\b(u)=\d_2(c)$ in $H^{k+1}(X, \BZ)$.
This fact follows from the commutativity of the following diagram.
\begin{equation*}
\begin{CD}
0 @>>> \BZ @>{n}>> \BZ@>>> \BZ_n@>>> 0\\
@. @V{=}VV @V{\frac 1n}VV @V{\frac 1n}VV @.\\
0 @>>> \BZ @>>> \BR@>>>  S^1 @>>> 0
\end{CD}
\end{equation*}

Since $H^k_{\text{$n$-torsion}}(X, S^1)=\hom(H_k(X,\BZ), \BZ_n)$, the
$n$-torsion elements in $\wh\BH^k(X)$ are given from the holonomy map point
of view by
$$
{\bold H}{\bold o}{\bold l}_{\text{$n$-torsion}}(X)\cong
\hom(H_k(X,\BZ), \BZ_n).
$$

\medskip
\noindent
{\bf The refined integer Stiefel-Whitney class.}
Suppose $F$ is an oriented rank $n$ real vector bundle.  Then the $k$th
Steifel-Whitney class of $F$ is a class $w_k(F)\in H^k(X;\BZ_2)$.
Suppose $k<n$ is even (otherwise twist by the orientation bundle of $F$).
The the image of $w_k(F)$ under the Bockstein map associated with $0\to
\BZ\overset{2}{\lra}\BZ\to \BZ_2\to0$ is the $(k+1)$st integer
Stiefel-Whitney class
$W_{k+1}(F)\in H^{k+1}_{\text{$2$-torsion}} (X;\BZ)$. This defines a unique
Stiefel-Whitney secondary class ${\widehat W}_{k+1}(F)\in H^{k}(X;S^1)
\subset\wh\BH^k(X)$, namely ${\widehat W}_{k+1}(F)$
 is the image of $w_k(F)$ under the mapping $H^k(X;\BZ_2)\to H^k(X; S^1)$
 as in the paragraph above.

\medskip
\noindent
{\bf Linear dependency currents and  Stiefel-Whitney sparks.}   Suppose $F$
is a rank $n$ real vector bundle.  Following [HZ$_1$], suppose
$\a=(\a_1,...,\a_{n-k})$ is an atomic collection of sections of $F$. Suppose
$F$ is oriented (or otherwise twist by the orientation bundle of $F$).
Assume $k$ is even. Then there exists a $d$-closed locally integrally flat
current $LD(\a)$ of degree $k+1$,  whose support is contained in the linear
dependency set of the collection of sections.  Furthermore, if $\b$ is
another atomic collection of sections, then
$$
LD(\a)-LD(\b)\ =\ dR
$$
where $R$ is a locally rectifiable current.

It is proven in [HZ$_1$] that the class of $LD(\a)$ in
$ H^{k+1} (X;\BZ)$ is $W_{k+1}(F)$, the  $(k+1)$st integer Stiefel-Whitney
class of $F$.

Now given a metric connection on $F$, Zweck [Z] constructs a canonical
spark $S_{k+1}(\a)$ which is a current with
$L^1_{\text{loc}}$-coefficients  satisfying the (flat) spark equation
$$
d S_{k+1}(\a)\ =\ -LD(\a)
$$
on $X$.  Furthermore, Zweck [Z]  proves a spark comparison formula
$$
 S_{k+1}(\a) -  S_{k+1}(\a')\ =\ T+dL
$$
with $T$ locally integrally flat and $L$ an
$L^1_{\text{loc}}$-current.  That is, there is a well defined  class
$$
[S_{k+1}(\a)] \in H^{k}(X;S^1) \subset \wh\BH^k(X)
$$
 Zweck [Z]  proves that this spark class, which refines $W_{k+1}(F)$
is the image of the $k$th Stiefel-Whitney class $w_k(F)\in H^k(X;\BZ_2)$
under the map   $(1/2):H^k(X;\BZ_2) \to H^k(X;S^1)$.  This proves that
the spark class of $S_{k+1}(\a)$ is the refined integer $(k+1)$st
Stiefel-Whitney class ${\widehat W}_{k+1}(F)$, discussed above.

\end{document}